\documentclass[12pt,a4paper,reqno]{amsart}  

\usepackage{mathtools}
\usepackage{url}
\usepackage[headings]{fullpage}
\usepackage{comment} 
\usepackage{datetime}  
\usepackage[english]{babel} 

\usepackage[full]{textcomp}
\usepackage[osf]{newtxtext} % osf for text, not math
\usepackage{cabin} % sans serif 
\usepackage[bigdelims,vvarbb]{newtxmath} % bb from STIX
\usepackage[cal=boondoxo]{mathalfa} % mathcal 
\usepackage{zlmtt}% latin modern typewriter

\usepackage{microtype} % migliora vari aspetti di microtipografia (giustificazione del testo ad esempio) 

\newcommand{\eps}{\varepsilon}

\newcommand{\dx}{\mathrm{d}} 
 
\newcommand{\I}{\mathcal{I}} 
\newcommand{\J}{\mathcal{J}}

\newcommand{\Etilde}{\widetilde{E}}

\newcommand{\Stilde}{\widetilde{S}}
 
\newcommand{\Vtilde}{\widetilde{V}}

\newcommand{\Odi}[1]{\Odip{}{#1}}

\newcommand{\Odim}[1]{\mathcal{O}\bigl(#1\bigr)}   
\newcommand{\Odip}[2]{\mathcal{O}_{#1}\left(#2\right)}
\newcommand{\Odipg}[2]{\mathcal{O}_{#1}\Bigl(#2\Bigr)}  
\newcommand{\Odipm}[2]{\mathcal{O}_{#1} (#2)}  
\newcommand{\odip}[2]{{o}_{#1}\left(#2\right)}
\newcommand{\odi}[1]{\odip{}{#1}}

\allowdisplaybreaks

\renewcommand{\qedsymbol}{$\square$}
\newenvironment{Proof}[1][Proof]{\par\noindent\textbf{#1.}~}
{\hfill\qedsymbol\smallskip\par}

\newtheoremstyle{sltheorems}% name
{10pt}%      Space above
{6pt}%      Space below
{\slshape}%         Body font
{}%         Indent amount (empty = no indent, \parindent = para indent)
{\bfseries}% Thm head font
{.}%        Punctuation after thm head
{.5em}%     Space after thm head: " " = normal interword space;
{\thmname{#1}\thmnumber{ #2}\thmnote{ (#3)}}
\theoremstyle{sltheorems} 
\newtheorem{Theorem}{Theorem}
\newtheorem{Lemma}{Lemma} 

\newtheoremstyle{remark}% name
{10pt}%      Space above
{6pt}%      Space below
{\rm} %         Body font
{}%         Indent amount (empty = no indent, \parindent = para indent)
{\bfseries}% Thm head font
{.}%        Punctuation after thm head
{.5em}%     Space after thm head: " " = normal interword space;
{\thmname{#1}\thmnumber{ #2}\thmnote{ (#3)}}
 \theoremstyle{remark}

\makeatletter

\def\env@Biggcases{%
  \let\@ifnextchar\new@ifnextchar
  \Biggl\lbrace
  \def\arraystretch{1.2}%
  \array{@{}l@{\quad}l@{}}%
}
\makeatother

\allowdisplaybreaks

\begin{document}

\title[one prime power and two squares of primes in short intervals]{Sums of one prime power and two squares of primes \\ in short intervals} 
\author[]{Alessandro Languasco \lowercase{and} Alessandro Zaccagnini}

%\date{\today, \currenttime}
%    \subjclass is required.
\subjclass[2010]{Primary 11P32; Secondary 11P55, 11P05}
\keywords{Waring-Goldbach problem, Laplace transforms}
\begin{abstract}
Let $k\ge 1$ be an integer.
We prove that a suitable
asymptotic formula for the average number  of representations
of  integers $n=p_{1}^{k}+p_{2}^{2}+p_{3}^{2}$,
where  $p_1,p_2,p_3$ are prime numbers,
holds in intervals shorter than the ones previously known.
\end{abstract}

\selectlanguage{english}
\maketitle

\section{Introduction}
The problem of representing an integer as a sum of a prime power and of two prime squares
is classical. 
It is conjectured that every sufficiently large $n$ subject to some congruence conditions
can be represented as $n= p_{1}^k +p_{2}^{2} + p_{3}^{2}$, where $k \ge 1$ is an integer. 
Let now $N$ be a large integer and denote by $E_k(N)$ the cardinality
of the set of integers not exceeding $N$ that satisfy the necessary
congruence conditions but can not be represented as the sum of a
$k$-th prime power and two prime squares.
Several results about $E_k(N)$ were obtained; the first one who proved
a non-trivial estimate
for $ E_k(N)$ was Hua \cite{Hua1938}. Later Schwarz
\cite{Schwarz1961a} and  several other authors
further  improved such an estimate;
we recall the contribution of 
 Leung-Liu \cite{LeungL1993},
Harman-Kumchev \cite{HarmanK2010},
L\"u \cite{Lu2006} and 
Li \cite{Li2012}.
Let $\eps>0$; so far the best known estimates are
$ E_1(N)\ll N^{1/3+\eps}$ 
by Zhao \cite{Zhao2014a},
$ E_2(N)\ll N^{17/20+\eps}$ 
by Harman-Kumchev \cite{HarmanK2010}, 
$ E_3(N)\ll N^{15/16+\eps}$ 
and, for $k\ge 4$,
$ E_k(N) \ll N^{1-1/(4k^2)+\eps}$  
both by Br\"udern \cite{Brudern2016}.
Let
\begin{equation}
\label{r-def}
  r_k(n)
  =
  \sum_{p_{1}^k + p_{2}^2 + p_{3}^{2} = n} \log p_{1} \log p_{2} \log p_{3}.
\end{equation} 

In this paper we study the average behaviour  of $r_k(n)$  
over  short intervals $[N,N+H]$, $H=\odi{N}$ thus generalising
our result in \cite{LanguascoZ2016a} which just deals with the case $k=1$.
  
\begin{Theorem}
\label{existence}  
 Let $N\ge 2$, $1\le H \le N$, $k\ge1$ be integers.
Then, for every $\eps>0$, there exists $C=C(\eps)>0$ such that
 \[
   \sum_{n = N+1}^{N + H}   r_k(n) = \frac{\pi }{4 }HN^{1/k}
   +
\Odipg{k}{HN^{1/k}  \exp \Bigl( -C \Bigl( \frac{ \log N}{\log  \log N} \Bigr)^{1/3} \Bigr)}
    \quad \textrm{as} \ N \to \infty,
 \] 
uniformly for  
$N^{1-5/(6k)+\eps} \le  H \le N^{1-\eps}$
for $k\ge 2$ and $N^{7/12+\eps} \le  H \le N^{1-\eps}$
for $k=1$.
% where $\Gamma$ is Euler's function. (???)
\end{Theorem}

It is worth remarking that 
the formula in Theorem \ref{existence} implies that
every interval $[N,N+H]$ contains an integer which is a sum of a prime
$k$-th power
and two prime squares, where 
$N^{1-5/(6k)+\eps} \le  H \le N^{1-\eps}$
for $k\ge 2$ and $N^{7/12+\eps} \le  H \le N^{1-\eps}$
for $k=1$.
In fact, for $k=1$ Zhao's estimate previously mentioned leads to better consequences than our Theorem \ref{existence}, but for $k\ge 2$ our result gives non-trivial information.

Assuming that the Riemann Hypothesis (RH) holds,
we prove that a suitable asymptotic formula for such an average of
$r_k(n)$ holds in much shorter intervals. We need the following
auxiliary function: let
\begin{equation}
\label{Ek-def}
E(k) : = 
\begin{cases}
N^{3/2} \log N + HN^{3/4} (\log N)^{3/2}& \text{if} \ k=1\\
N \log N+ HN^{1/4} (\log N)^{2} & \text{if} \ k=2\\
N^{5/6} \log N+ HN^{1/4} \log N  + H^{1/2} N^{1/2} \log N & \text{if} \ k=3\\
N^{3/4+1/k} \log N & \text{if} \ k\ge 4.
\end{cases}
\end{equation}
We have the following
\begin{Theorem}
\label{existence-RH}
Assume the Riemann Hypothesis (RH).
 Let $N\ge 2$, $1\le H \le N$, $k\ge 2$ be integers. 
 We have
 \[
   \sum_{n = N+1}^{N + H}   r_k(n) = \frac{\pi  }{4  }HN^{1/k}
   +
\Odipm{k}{H^2 N^{1/k-1}+H^{1/2}N^{1/2+1/(2k)}(\log N)^2 + E(k)}
 \] 
as $N \to \infty$,
uniformly for  $\infty(N^{1-1/k}(\log N)^4)\le H \le \odi{N}$, where $f=\infty(g)$ means $g=\odi{f}$,  $\Gamma$ is Euler's function and $E(k)$ is defined in \eqref{Ek-def}. 
\end{Theorem}
We remark that a  version for  $k=1$ of Theorem  \ref{existence-RH}  was obtained in \cite{LanguascoZ2016a}.
We further remark that 
the formula in Theorem \ref{existence-RH} implies that
every interval $[N,N+H]$ contains an integer which is a sum of a prime power
and two prime squares, where $\infty(N^{1-1/k}(\log N)^4)\le   H =\odi{N}$.

%For $k=2,3$ we can remove the condition $H\ge \infty(N^{3/4}\log N)$ by replacing
%$r_k(n)$ with
%$
%  R_k(n)
%  =
%  \sum_{m_{1}^k + m_{2}^2 + m_{3}^{2} = n}  \Lambda(m_{1})  \Lambda(m_{2}) \Lambda(m_{3})
%  $, where $\Lambda$ is the von Mangoldt function, 
%  since in this way we can discard the contribution 
%  of $\I_{4}, \I_{5}$ of sections \ref{sect:I5-estim-RH}-\ref{sect:I4-estim-RH} below.
  
The proofs of both Theorems \ref{existence}-\ref{existence-RH} use the original
Hardy-Littlewood settings of the circle method to exploit  
the easier main term treatment they allow (comparing
with the one which would follow
using  Lemmas 2.3 and 2.9 of Vaughan \cite{Vaughan1997}).

It is worth remarking that the expected best result using circle method techniques
is $H \ge  N^{1-1/k}$; so our Theorem \ref{existence-RH}, under the assumption
of the Riemann  Hypothesis, 
comes very close to this bound.
We also obtained similar results in \cite{LanguascoZ2017b} and
\cite{LanguascoZ2017c}.

\section{Notation and Lemmas}
 Let   $e(\alpha) = e^{2\pi i\alpha}$,
$\alpha \in [-1/2,1/2]$, $L=\log N$, $z= 1/N-2\pi i\alpha$,
\begin{equation*} 
%\label{tilde-main-defs}
\Stilde_\ell(\alpha) = \sum_{n=1}^{\infty} \Lambda(n) e^{-n^{\ell}/N} e(n^{\ell}\alpha)   \quad
\textrm{and}
\quad
\Vtilde_\ell(\alpha) = \sum_{p=2}^{\infty} \log p \, e^{-p^{\ell}/N} e(p^{\ell}\alpha).
\end{equation*} 
We remark that 
\begin{equation}
\label{z-estim}
\vert z\vert^{-1} \ll \min \bigl(N, \vert \alpha \vert^{-1}\bigr).
\end{equation}
 We further set
\begin{equation}
\notag 
   U(\alpha,H)
  = 
  \sum_{m=1}^H e(m \alpha) ,
\end{equation}
and, moreover, we also have the usual 
numerically explicit inequality
\begin{equation}
\label{UH-estim}
\vert U(\alpha,H) \vert
\le
\min \bigl(H; \vert \alpha \vert^{-1}\bigr),
\end{equation}
see, \emph{e.g.}, on page 39 of Montgomery \cite{Montgomery1994}. %
We list now the needed preliminary results.
 \begin{Lemma}[Lemma 3 of  \cite{LanguascoZ2016b}]
\label{tilde-trivial-lemma}
Let $\ell\ge 1$ be an integer. Then
\(
\vert \Stilde_{\ell}(\alpha)- \Vtilde_{\ell}(\alpha) \vert 
\ll_{\ell}
 N^{1/(2\ell)}  .
\)
\end{Lemma}

\begin{Lemma}
\label{Linnik-lemma}
Let $\ell \ge 1$ be an integer, $N \ge 2$  and $\alpha\in [-1/2,1/2]$.
Then
\begin{equation*}
%\label{expl-form}
\Stilde_{\ell}(\alpha)  
= 
\frac{\Gamma(1/\ell)}{\ell z^{1/\ell}}
- 
\frac{1}{\ell}\sum_{\rho}z^{-\rho/\ell}\Gamma (\rho/\ell ) 
+
\Odip{\ell}{1},
\end{equation*}
where $\rho=\beta+i\gamma$ runs over
the non-trivial zeros of $\zeta(s)$.
\end{Lemma}
\begin{Proof}
It follows the line of  Lemma 2 of  \cite{LanguascoZ2016a}; we just 
correct an oversight in its  proof. In eq. (5) on page 48 of 
\cite{LanguascoZ2016a} the term 
\(
-  
  \sum_{m=1}^{\ell \sqrt{3}/4} \Gamma (- 2m/\ell ) z^{2m/\ell}
\)
is missing. Its estimate is trivially $\ll_{\ell} \vert z \vert^{\sqrt{3}/2} \ll_{\ell} 1$.
Hence such an oversight does not affect the final result of  
Lemma 2 of  \cite{LanguascoZ2016a}.
\end{Proof}

\begin{Lemma} [Lemma 4 of \cite{LanguascoZ2016a}]
 \label{Laplace-formula}
Let $N$ be a positive integer and 
$\mu > 0$.
Then
\[
  \int_{-1 / 2}^{1 / 2} z^{-\mu} e(-n \alpha) \, \dx \alpha
  =
  e^{- n / N} \frac{n^{\mu - 1}}{\Gamma(\mu)}
  +
  \Odipg{\mu}{\frac{1}{n}},
\]
uniformly for $n \ge 1$.
\end{Lemma}

\begin{Lemma}
 \label{LP-Lemma-gen} 
Let $\eps$ be an arbitrarily small
positive constant,  $\ell \ge 1$ be an integer, $N$ be a
sufficiently large integer and $L= \log N$. Then there exists a positive constant 
$c_1 = c_{1}(\eps)$, which does not depend on $\ell$, such that 
\[
\int_{-\xi}^{\xi} \,
\Bigl\vert
\Stilde_\ell(\alpha) - \frac{\Gamma(1/\ell)}{\ell z^{1/\ell}}
\Bigr\vert^{2}
\dx \alpha 
\ll_{\ell}
 N^{2/\ell-1} \exp \Big( - c_{1}  \Big( \frac{L}{\log L} \Big)^{1/3} \Big) 
\]
uniformly for $ 0\le \xi < N^{-1 +5/(6\ell) - \eps}$.
Assuming RH we get 
\[
\int_{-\xi}^{\xi} \,
\Bigl\vert
\Stilde_\ell(\alpha) - \frac{\Gamma(1/\ell)}{\ell z^{1/\ell}}
\Bigr\vert^{2}
\dx \alpha 
\ll_{\ell}
N^{1/\ell}\xi L^{2}
\]
uniformly  for  $0 \le \xi \le 1/2$.
\end{Lemma} 
\begin{Proof}
It follows the line of Lemma 3 of \cite{LanguascoZ2016a} and 
Lemma 1 of \cite{LanguascoZ2016b}; we just 
correct an oversight in their  proofs. Both eq. (8) on page 49 of 
\cite{LanguascoZ2016a}  and eq. (6) on page 423
of \cite{LanguascoZ2016b} should read as 
\[
\int_{1/N}^{\xi}
\Big \vert\sum_{\rho\colon \gamma > 0}z^{-\rho/\ell}\Gamma (\rho/\ell ) \Big \vert^2 
\dx \alpha 
\le
\sum_{k=1}^K
\int_\eta^{2\eta} \Big \vert\sum_{\rho\colon \gamma > 0}z^{-\rho/\ell}\Gamma (\rho/\ell ) \Big \vert^2 \dx \alpha,
\]
where $\eta=\eta_k= \xi/2^k$, $1/N\le \eta \le \xi/2$  and $K$ is a suitable integer satisfying $K=\Odi{L}$. 
The remaining part of the proofs are left untouched. 
Hence such  oversights do not affect the final result of  
Lemma 3 of \cite{LanguascoZ2016a} and 
Lemma 1 of \cite{LanguascoZ2016b}.
\end{Proof}

\begin{Lemma}[Lemma 2 of \cite{LanguascoZ2016b}]
\label{zac-lemma-series}
Let $\ell\ge 2 $ be an integer and $0<\xi\le 1/2$. Then
\[ 
\int_{-\xi}^{\xi} 
|\Stilde_{\ell}(\alpha)|^2 \ \dx\alpha 
\ll_{\ell}
\xi N^{1/\ell} L  +
\begin{cases}
L^{2} & \text{if}\ \ell =2\\
1 & \text{if}\ \ell > 2
\end{cases}
\
\text{and}
\
\int_{-\xi}^{\xi} 
|\Vtilde_{\ell}(\alpha)|^2 \ \dx\alpha 
\ll_{\ell}
\xi N^{1/\ell} L  +
\begin{cases}
L^{2} & \text{if}\ \ell =2\\
1 & \text{if}\ \ell > 2.
\end{cases}
\]
\end{Lemma}  
\begin{Proof}
The first   part was proved in Lemma 2 of \cite{LanguascoZ2016b}.
For the second part we argue analogously.
We use Corollary 2 of Montgomery-Vaughan  \cite{MontgomeryV1974}  
with $T=\xi$, 
 $a_r=\log(r) \exp(-r^\ell/N)$ if $r$ is prime,  $a_r= 0$ otherwise
 and $\lambda_r= 2\pi r^\ell$.
By the Prime Number Theorem we get
\begin{align*}
\int_{0}^{\xi} \vert \Vtilde_{\ell}(\alpha) \vert^2\, \dx \alpha 
&=
\sum_{p} \log^2 (p)  e^{-2p^\ell/N} \bigl(\xi +\Odim{\delta_p^{-1}}\bigr)
\ll_\ell
\xi N^{1/\ell} L
+
\sum_{p} \log^2 (p) p^{1-\ell} e^{-2p^\ell/N}
\end{align*}
since $\delta_{r} = \lambda_r - \lambda_{r-1} \gg_{\ell} r^{\ell-1}$. 
The last  term is $\ll_{\ell}1$ if $\ell >2$ and $\ll L^{2}$ otherwise.
The second part of Lemma \ref{zac-lemma-series} follows.
\end{Proof}  
%
\begin{comment}
In the unconditional case a crucial role is played by 
\begin{Lemma}[Hua]
\label{Hua-lemma-series}
%
There exists a suitable absolute constant $A>0$ such that
%
\[ 
\int_{-1/2}^{1/2} 
|\Stilde_{3}(\alpha)|^8 \ \dx\alpha 
\ll
 N^{5/3} L^A.
 \]
%
\end{Lemma} 
\begin{Proof}
The result follows  via partial integration from the analogue one for the
$S_{3}(\alpha;P) = \sum_{n\le P} \Lambda(n) e(n^3\alpha)$
in the form described
in section 4.3 of Deshouillers \cite{Deshouillers1990}.
\end{Proof}
\end{comment}

In the following we will also need a fourth-power average of $\Stilde_2(\alpha)$.

\begin{Lemma}[Lemma 6 of  \cite{LanguascoZ2016a}]
\label{Hua-Rieger-Lemma}
We have 
\[
\int_{-1/2}^{1/2} \,
\vert
\Stilde_2(\alpha) 
\vert^{4}
\, \dx \alpha 
\ll
N L^{2}.
\]
\end{Lemma} 

\section{Proof of Theorem \ref{existence}}
 \label{unconditional}
Let $\eps>0$ and  $H>2B$, where
\begin{equation}
\label{B-def} 
B= B(N,d)=  \exp \Bigl( d   \Bigl( \frac{\log N}{\log \log N} \Bigr)^{1/3} \Bigr),
\end{equation} 
where   $d=d(\eps)>0$ will be chosen later. 
Recalling \eqref{r-def}, we may write
\begin{align*}  
  \sum_{n = N+1}^{N + H} 
  e^{-n / N} r_k(n)
  =
  \int_{-1/2}^{1/2}
    \Vtilde_k(\alpha) \Vtilde_2(\alpha)^{2}  U(-\alpha,H) e(-N\alpha) \, \dx \alpha.
\end{align*} 
We find it also convenient to set
\begin{equation}
\label{def-Etilde}
  \Etilde_{\ell}(\alpha)
  : =
  \Stilde_\ell(\alpha) - \frac{\Gamma(1/\ell)}{\ell z^{1/\ell}}.
\end{equation}

Letting $I(B,H):=[-1/2,-B/H]\cup  [B/H, 1/2]$,
using Lemma \ref{Linnik-lemma}, recalling that
 $\Gamma(1/2)=\pi^{1/2}$, we can write 
\begin{align}
\notag
  \sum_{n = N+1}^{N + H} 
  e^{-n / N} r_k(n)
  &= 
    \int_{-B/H}^{B/H}
    \frac{\pi\Gamma(1/k)}{4k z^{1+1/k}} 
      U(-\alpha,H) e(-N\alpha) \, \dx \alpha
\\
\notag
&\hskip1cm +
       \int_{-B/H}^{B/H}
      \frac{\Gamma(1/k)}{k z^{1/k}} \Bigl( \Stilde_2(\alpha)^{2}-\frac{\pi}{4z}\Bigr) 
    U(-\alpha,H) e(-N\alpha) \, \dx \alpha
\\
\notag
&\hskip1cm +
\int_{-B/H}^{B/H}
    \Etilde_k(\alpha)  \Stilde_2(\alpha)^{2} U(-\alpha,H) e(-N\alpha) \, \dx \alpha \\
\notag
&\hskip1cm +
  \int_{-1/2}^{1/2}
   \Vtilde_k(\alpha) ( \Vtilde_2(\alpha)^{2} - \Stilde_2(\alpha)^{2}) U(-\alpha,H) e(-N\alpha) \, \dx \alpha 
\\\notag
&\hskip1cm +
  \int_{-1/2}^{1/2}
   \Stilde_2(\alpha)^2 ( \Vtilde_k(\alpha) - \Stilde_k(\alpha) ) U(-\alpha,H) e(-N\alpha) \, \dx \alpha 
   \\\notag
&\hskip1cm +
     \int_{I(B,H)}
    \Stilde_k(\alpha) \Stilde_2(\alpha)^{2}  U(-\alpha,H) e(-N\alpha) \, \dx \alpha
    \\
    \label{main-split}
& =
    \J_{1}+\J_{2}+\J_{3}+ \J_{4}+ \J_{5}+ \J_{6},
\end{align}
say. Now we evaluate these terms.

\subsection{Evaluation of $\J_{1}$}
\label{sect:J1-eval}
Using Lemma~\ref{Laplace-formula} and \eqref{z-estim} we immediately get 
\begin{align}
\notag
\J_{1}
& =
\frac{\pi \Gamma(1/k)}{4 k\Gamma(1+1/k)}  \sum_{n = N+1}^{N + H}  n^{1/k}e^{-n/N} 
+ \Odipg{k}{\frac{H}{N}} 
+ \Odipg{k}{\int_{B/H}^{1/2} \frac{\dx \alpha}{\alpha^{2+{1/k}}}}
\\&
\label{J1-eval}
=
\frac{\pi }{4   e} HN^{1/k} +
\Odipg{k}{H^2 N^{1/k-1}+ N^{1/k} +\Bigl(\frac{H}{B}\Bigr)^{1+1/k}}.
\end{align}

\subsection{Estimation of $\J_{6}$}
\label{sect:J6-estim}
 
Using $\Stilde_k(\alpha) \ll_k N^{1/k}$,  
\eqref{UH-estim}, Lemma \ref{zac-lemma-series}
and a partial integration, we obtain 
that
\begin{align}
\notag
\J_6 
&\ll_k
N^{1/k} 
\int_{B/H}^{1/2}
\frac{\vert \Stilde_2(\alpha)\vert ^{2}}{\alpha}   \, \dx \alpha
\\ \notag
& \ll_k
N^{1/k} 
\Bigl[
\frac{H}{B} \Bigl(\frac{N^{1/2}B L}{H} + L^2 \Bigr) 
+ 
N^{1/2}L 
+
\int_{B/H}^{1/2} \frac{N^{1/2}\xi L+L^2}{\xi^2}\ \dx \xi
\Bigr]
\\
\label{J6-estim}
&
\ll_k
N^{1/k} L^2 \Bigl(N^{1/2} + \frac{H}{B}\Bigr)
\end{align}
which, comparing with \eqref{J1-eval}, is under control for $H=\infty(N^{1/2}L^2)$
and $B = \infty(L^2)$ (which is fine thanks to \eqref{B-def}).

\subsection{Estimation of $\J_{5}$}
\label{sect:J5-estim}
By Lemmas \ref{tilde-trivial-lemma} and \ref{zac-lemma-series},
\eqref{UH-estim} and a partial integration
  we get  
\begin{align}
\notag
\J_{5}
&\ll
  \int_{-1/2}^{1/2}
   \vert \Stilde_2(\alpha)\vert ^2  \vert \Vtilde_k(\alpha) - \Stilde_k(\alpha) \vert \vert U(-\alpha,H) \vert\, \dx \alpha 
\\&
\notag
\ll_k
HN^{1/(2k)} 
\int_{-1/H}^{1/H}  
\vert \Stilde_2(\alpha)\vert ^2   \ \dx \alpha 
+
N^{1/(2k)} 
\int_{1/H}^{1/2}  
\frac{\vert \Stilde_2(\alpha)\vert ^2}{\alpha}\ \dx \alpha 
\\
\notag
&
\ll_k
HN^{1/(2k)} \Bigl(\frac{N^{1/2}L}{H}+ L^2 \Bigr)
+
N^{1/(2k)} 
\Bigl( 
N^{1/2}L+ HL^2  
+
\int_{1/H}^{1/2}  
\frac{N^{1/2} \xi L +L^2}{\xi^2}\ \dx \xi 
 \Bigr)
 \\
 \label{J5-estim-series}
&
\ll_k
N^{1/(2k)} ( N^{1/2} + H ) L^2.
\end{align}
which, comparing with \eqref{J1-eval}, is under control for $H=\infty(N^{1/2-1/(2k)}L^2)$.

\subsection{Estimation of $\J_{4}$}
\label{sect:J4-estim}

Using the identity $f^2-g^2= 2f(f-g)- (f-g)^2$, Lemma \ref{tilde-trivial-lemma} 
and $\Vtilde_k(\alpha)  \ll_k N^{1/k}$,
we  have
\begin{align}
\notag
\Vtilde_k(\alpha) ( \Vtilde_2(\alpha)^{2} - \Stilde_2(\alpha)^{2}) 
& \ll_k
\vert \Vtilde_k(\alpha) \vert 
\bigl(  
\vert \Vtilde_2(\alpha) \vert   
\vert \Vtilde_{2}(\alpha)  -\Stilde_{2}(\alpha) \vert
+
 \vert \Vtilde_{2}(\alpha)  -\Stilde_{2}(\alpha) \vert^2
\bigr)
\\
&
\notag
\ll_k 
N^{1/4}
\vert \Vtilde_k(\alpha) \vert   
\vert \Vtilde_2(\alpha) \vert    
+
N^{1/2+1/k}.
\end{align}
Clearly we have
\begin{align}
\notag
\J_4 
&\ll_k 
N^{1/4}
\int_{-1/2}^{1/2} 
\vert \Vtilde_k(\alpha) \vert   
\vert \Vtilde_2(\alpha) \vert   
\vert U(-\alpha,H)
\vert \, \dx \alpha
+
N^{1/2+1/k}
\int_{-1/2}^{1/2}  
\vert U(-\alpha,H)
\vert \, \dx \alpha
\\&
\label{J4-split}
= K_1+K_2,
\end{align}
say.
Using \eqref{UH-estim} we get
\begin{equation}
\label{UH-average}
\int_{-1/2}^{1/2} 
\vert U(-\alpha,H)
\vert \, \dx \alpha
\ll
\int_{-1/H}^{1/H}  
H  \ \dx \alpha 
+
\int_{1/H}^{1/2}  
\frac{\dx \alpha}{\alpha}
\ll 
L
\end{equation}
and hence, by \eqref{J4-split}-\eqref{UH-average}, we can write
\begin{equation}
\label{K2-estim}
K_2 
\ll_k
N^{1/2+1/k} L,
\end{equation}
for every $k\ge 1$.

Now we estimate $K_1$; depending on $k$, we need to perform
different computations.

Let $k=1$.
Using the Cauchy-Schwarz inequality, the Prime Number Theorem, 
 \eqref{UH-estim} and Lemma \ref{zac-lemma-series},
 we obtain that
\begin{align}
\notag  
 K_1  
&\ll
N^{1/4}
\Bigl(
 \int_{-1/2}^{1/2} 
\vert  \Vtilde_1 (\alpha) \vert^2 \, \dx \alpha 
\Bigr)^{1/2}
\Bigl(
 \int_{-1/2}^{1/2} 
\vert  \Vtilde_2 (\alpha) \vert^2
\vert U(-\alpha,H)
\vert^2 \, \dx \alpha 
\Bigr)^{1/2}
\\
\notag
&
\ll
N^{3/4}L^{1/2}
\Bigl[
H^2 \Bigl(\frac{N^{1/2} L}{H} + L^2 \Bigr) 
+ 
N^{1/2}L 
+
\int_{1/H}^{1/2} \frac{N^{1/2}\xi L+L^2}{\xi^3}\ \dx \xi
\Bigr] ^{1/2}
\\
\label{K1-estim-k=1}
&
\ll
H N^{3/4}L^{3/2} + 
H^{1/2} N L.
\end{align}
Hence, by \eqref{J4-split} and \eqref{K2-estim}-\eqref{K1-estim-k=1}, for $k=1$ we get 
\begin{equation}
\label{J4-estim-series-k=1}
\J_{4}
\ll
N^{3/2} L +  
H N^{3/4}L^{3/2}.
\end{equation}

Let $k=2$.
Using 
 \eqref{UH-estim} and Lemma \ref{zac-lemma-series},
 we obtain that
\begin{align}
\notag 
 K_1 & 
\ll
N^{1/4}  
 \int_{-1/2}^{1/2} 
\vert  \Vtilde_2 (\alpha) \vert^2
\vert U(-\alpha,H)
\vert \, \dx \alpha  
\\
\notag
&
\ll
N^{1/4}
\Bigl[
H \Bigl(\frac{N^{1/2} L}{H} + L^2 \Bigr) 
+ 
N^{1/2}L 
+
\int_{1/H}^{1/2} \frac{N^{1/2}\xi L+L^2}{\xi^2}\ \dx \xi
\Bigr] 
\\
\label{K1-estim-k=2}
&
\ll
H N^{1/4}L^{2} + 
 N^{3/4} L^2.
\end{align}
Hence, by \eqref{J4-split}, \eqref{K2-estim} and \eqref{K1-estim-k=2}, for $k=2$ we get 
\begin{equation}
\label{J4-estim-series-k=2}
\J_{4}
\ll
N L +   H N^{1/4}L^{2} .
\end{equation}

Let $k=3$.
Using the Cauchy-Schwarz estimate, \eqref{UH-estim} and Lemma \ref{zac-lemma-series},
 we obtain that
\begin{align}
\notag  
 K_1  &   
\ll
N^{1/4}  
\Bigl(
 \int_{-1/2}^{1/2} 
\vert  \Vtilde_3 (\alpha) \vert^2
\vert U(-\alpha,H)
\vert \, \dx \alpha 
\Bigr)^{1/2}
\Bigl(
 \int_{-1/2}^{1/2} 
\vert  \Vtilde_2 (\alpha) \vert^2
\vert U(-\alpha,H)
\vert \, \dx \alpha 
\Bigr)^{1/2}
\\
\notag
&
\ll
N^{1/4}
\Bigl[
H \Bigl(\frac{N^{1/3} L}{H} + 1 \Bigr) 
+ 
N^{1/3}L 
+
\int_{1/H}^{1/2} \frac{N^{1/3}\xi L+1}{\xi^2}\ \dx \xi
\Bigr]^{1/2}
\\
\notag
&
\hskip2cm
\times
\Bigl[
H \Bigl(\frac{N^{1/2} L}{H} + L^2 \Bigr) 
+ 
N^{1/2}L 
+
\int_{1/H}^{1/2} \frac{N^{1/2}\xi L+L^2}{\xi^2}\ \dx \xi
\Bigr]^{1/2}
\\
\label{K1-estim-k=3}
&
\ll
H N^{1/4}L  + H^{1/2} N^{1/2}L +
 N^{2/3} L^2.
\end{align}
Hence, by  \eqref{J4-split}, \eqref{K2-estim} and \eqref{K1-estim-k=3}, for $k=3$ we get 
\begin{equation}
\label{J4-estim-series-k=3}
\J_{4}
\ll
N^{5/6} L +   H N^{1/4}L + H^{1/2} N^{1/2}L  .
\end{equation}

Let $k\ge 4$.
By \eqref{J4-split}, $\Vtilde_k(\alpha)  \ll_k N^{1/k}$ and \eqref{UH-average} we can  write
that 
\begin{equation}
 \label{J4-estim-series-k>=4}
\J_{4}
\ll_k
N^{3/4+1/k}
\int_{-1/2}^{1/2} 
\vert U(-\alpha,H)
\vert \, \dx \alpha 
\ll_k
N^{3/4+1/k} L.
\end{equation}

Summing up, from \eqref{J4-estim-series-k=1}, \eqref{J4-estim-series-k=2}
and
\eqref{J4-estim-series-k=3}-\eqref{J4-estim-series-k>=4} we can write that
\begin{equation}
 \label{J4-estim-series}
\J_{4} \ll_k E(k),
\end{equation}
where $E(k)$ is defined in \eqref{Ek-def}.

\subsection{Estimation of $\J_{2}$}
\label{sect:I2-estim}
Now we estimate $\J_{2}$.
Using the identity $f^{2}-g^{2} = 2f(f-g) -(f-g)^{2}$ we obtain 
\begin{equation}
\label{J2-estim1}
\J_{2} 
\ll_k
 \int_{-B/H}^{B/H}
      \vert \Etilde_{2}(\alpha)  \vert
        \frac{\vert U(\alpha,H)\vert}{\vert z\vert^{1/2+1/k}}
     \, \dx \alpha
+
 \int_{-B/H}^{B/H}
       \vert \Etilde_{2}(\alpha)  \vert^{2}  
       \frac{\vert U(\alpha,H) \vert}{\vert z\vert^{1/k}} 
        \, \dx \alpha 
    =
   I_{1} +I_{2},
\end{equation}
say.
Using   \eqref{z-estim},  \eqref{UH-estim}, Lemma \ref{LP-Lemma-gen} and a partial integration argument  we obtain   that,  for every $\eps>0$, there exists $c_1=c_1(\eps)>0$ such that
\begin{equation} 
 \label{I2-estim}
I_{2} 
\ll_k
HN^{1/k}
 \int_{-B/H}^{B/H} \vert \Etilde_{2}(\alpha)  \vert^{2}   \, \dx \alpha  
\ll_k
 HN^{1/k}
\exp \Big( -  c_{1}  \Big( \frac{L}{\log L} \Big)^{1/3} \Big)
\end{equation} 
provided that $H\ge B N^{7/12+\eps}$.
Using  the Cauchy-Schwarz inequality and arguing as for $I_{2}$ we get 
\begin{equation}
  \label{I1-estim}  
I_{1}
\ll_k
H   
\Bigl(
 \int_{-B/H}^{B/H}  \vert \Etilde_{2}(\alpha)  \vert^{2}   \, \dx \alpha 
 \Bigr)^{1/2} 
  \Bigl(   \int_{-B/H}^{B/H}  \frac{\dx \alpha}{\alpha^{1+2/k}}\Bigr)^{1/2}  
 \ll_k
H   N^{1/k} \exp \Big( - \frac{c_{1}}{2}  \Big( \frac{L}{\log L} \Big)^{1/3} \Big),
\end{equation}
provided that $H\ge B N^{7/12+\eps}$.
Inserting \eqref{I2-estim}-\eqref{I1-estim} into \eqref{J2-estim1} we finally obtain
\begin{equation}
\label{J2-estim-final}
\J_{2} \ll_k H   N^{1/k} \exp \Big( - \frac{c_{1}}{2}  \Big( \frac{L}{\log L} \Big)^{1/3} \Big),
\end{equation} 
provided that $H\ge B N^{7/12+\eps}$.

\subsection{Estimation of $\J_{3}$}
\label{sect:I3-estim}
Now we estimate $\J_{3}$.
By the Cauchy-Schwarz  inequality,  \eqref{UH-estim}, Lemmas
\ref{Hua-Rieger-Lemma} and Lemma \ref{LP-Lemma-gen},
 we obtain   that,  for every $\eps>0$, there exists $c_1=c_1(\eps)>0$ such that
\begin{align}
\notag
\J_{3}
&\ll_k
\Bigl( \int_{-1/2}^{1/2}  \vert\Stilde_2(\alpha)\vert^{4}   \, \dx \alpha \Bigr)^{1/2}
\Bigl( 
\int_{-B/H}^{B/H}  \vert \Etilde_{k}(\alpha)  \vert^{2} \vert U(\alpha,H)\vert^2 
 \, \dx \alpha 
 \Bigr)^{1/2}
\\
 \label{J3-estim-final}
&\ll_k
HN^{1/2}L 
\Bigl( 
\int_{-B/H}^{B/H}  \vert \Etilde_{k}(\alpha)  \vert^{2} \, \dx \alpha 
 \Bigr)^{1/2}
%\ll_k
%HN^{1/2}L   N^{1/k-1/2} \exp \Big( - \frac{c_{1}}{2}  \Big( \frac{L}{\log L} \Big)^{1/3} \Big) 
%\\
\ll_k
H N^{1/k} \exp \Big( - \frac{c_{1}}{2}  \Big( \frac{L}{\log L} \Big)^{1/3} \Big) ,
 \end{align}
provided that $H\ge B N^{1-5/(6k)+\eps}$.

\subsection{Final words}
Let $k\ge 2$.
By \eqref{main-split}-\eqref{J5-estim-series}, \eqref{J4-estim-series} and \eqref{J2-estim-final}-\eqref{J3-estim-final}
we have  that,  for every $\eps>0$, there exists $c_1=c_1(\eps)>0$ such that
\begin{align}
\notag
  \sum_{n = N+1}^{N + H} 
  e^{-n / N} r_k(n)
   &= \frac{\pi }{4   e} HN^{1/k}
  \\
  \label{almost-done}
&+   
\Odipg{k}{ H N^{1/k} \exp \Big( - \frac{c_{1}}{2}  \Big( \frac{L}{\log L} \Big)^{1/3} \Big) +\frac{H}{B} N^{1/k} L^2 + E(k)}
\end{align}
provided that $H\ge B N^{1-5/(6k)+\eps}$.
The second error term is dominated by the first one by choosing
$d=c_1$ in \eqref{B-def}. 
So from now on we have $H \ge N^{1-5/(6k)+\eps}$ for $k\ge 2$.
The third error term in \eqref{almost-done} is now dominated by the first. 

Let $k=1$. In this case 
\eqref{almost-done} holds  provided that 
$H\ge B N^{7/12+\eps}$ and
the second error term is dominated by the first one by choosing
$d=c_1$ in \eqref{B-def}. Hence, for $k=1$, we get
that $H\ge N^{7/12+\eps}$. 
The third error term in \eqref{almost-done} is now dominated by the first. 

Summing up,  for every $k\ge 1$
we can write  that,  for every $\eps>0$, there exists $C=C(\eps)>0$ such that
\begin{equation*}
%\label{almost-done-2}
  \sum_{n = N+1}^{N + H} 
  e^{-n / N} r_k(n)
  = \frac{\pi }{4  e} HN^{1/k}
+\Odipg{k}{ H N^{1/k} \exp \Big( - C  \Big( \frac{L}{\log L} \Big)^{1/3} \Big)}
\end{equation*}
provided that $H \ge N^{1-5/(6k)+\eps}$ for $k\ge 2$ and
 $H\ge N^{7/12+\eps}$ for $k=1$.
{}From  $e^{-n/N}=e^{-1}+ \Odi{H/N}$ for $n\in[N+1,N+H]$, $1\le H \le N$,
we get  that,  for every $\eps>0$, there exists $C=C(\eps)>0$ such that
\begin{equation*}
   \sum_{n = N+1}^{N + H} 
r_k(n)
=
 \frac{\pi }{4  } HN^{1/k}
+\Odipg{k}{ H N^{1/k} \exp \Big( - C\Big( \frac{L}{\log L} \Big)^{1/3} \Big)}
\\
+
  \Odipg{k}{\frac{H}{N}\sum_{n = N+1}^{N + H} r_k(n)
}
\end{equation*}
provided that $H \le N$, $H \ge N^{1-5/(6k)+\eps}$ for $k\ge 2$ and
 $H\ge N^{7/12+\eps}$ for $k=1$.
Using $e^{n/N}\le  e^{2}$ 
and \eqref{almost-done},
 the last error term is
$\ll_k H^2N^{1/k-1}$.
Hence we get
\begin{equation*}
%\label{th-series}
   \sum_{n = N+1}^{N + H} 
r_k(n)
=
 \frac{\pi }{4 } HN^{1/k}
+
\Odipg{k}{ H N^{1/k} \exp \Big( - C  \Big( \frac{L}{\log L} \Big)^{1/3} \Big)}
\end{equation*}
uniformly for  $N^{1-5/(6k)+\eps} \le  H \le N^{1-\eps}$
if $k\ge 2$ and for $N^{7/12+\eps} \le  H \le N^{1-\eps}$
if $k=1$.
Theorem \ref{existence} follows.

\section{Proof of Theorem \ref{existence-RH}}

Let $k\ge 2$, $H\ge 2$, $H=\odi{N}$  be an  integer. 
We recall that we set $L = \log N$ for brevity. 
From now on we assume that RH holds.
we may write
\begin{align*}  
  \sum_{n = N+1}^{N + H} 
  e^{-n / N} r_k(n)
  =
  \int_{-1/2}^{1/2}
    \Vtilde_k(\alpha) \Vtilde_2(\alpha)^{2}  U(-\alpha,H) e(-N\alpha) \, \dx \alpha.
\end{align*} 

In this conditional case we can simplify the setting.
Using Lemma \ref{Linnik-lemma},  
 recalling definition \eqref{def-Etilde} and that
 $\Gamma(1/2)=\pi^{1/2}$,
%\(
%\Etilde_{\ell}(\alpha) : =\Stilde_\ell(\alpha) - \frac{\Gamma(1/\ell)}{\ell z^{1/\ell}}
%\),
we can write 
\begin{align}
\notag
  \sum_{n = N+1}^{N + H} 
  e^{-n / N} r_k(n)
  &= 
    \int_{-1/2}^{1/2}
    \frac{\pi\Gamma(1/k)}{4k z^{1+1/k}} 
      U(-\alpha,H) e(-N\alpha) \, \dx \alpha
\\
\notag
& \hskip1cm +
       \int_{-1/2}^{1/2}
      \frac{\Gamma(1/k)}{k z^{1/k}} \Bigl( \Stilde_2(\alpha)^{2}-\frac{\pi}{4z}\Bigr) 
    U(-\alpha,H) e(-N\alpha) \, \dx \alpha
\\
\notag
& \hskip1cm +
\int_{-1/2}^{1/2}
    \Etilde_k(\alpha)  \Stilde_2(\alpha)^{2} U(-\alpha,H) e(-N\alpha) \, \dx \alpha \\
\notag
& \hskip1cm +
  \int_{-1/2}^{1/2}
   \Vtilde_k(\alpha) ( \Vtilde_2(\alpha)^{2} - \Stilde_2(\alpha)^{2}) U(-\alpha,H) e(-N\alpha) \, \dx \alpha 
\\\notag
& \hskip1cm +
  \int_{-1/2}^{1/2}
   \Stilde_2(\alpha)^2 ( \Vtilde_k(\alpha) - \Stilde_k(\alpha) ) U(-\alpha,H) e(-N\alpha) \, \dx \alpha 
\\
    \label{main-split-RH}
& =
    \I_{1}+\I_{2}+\I_{3}+ \I_{4}+ \I_{5},
\end{align}
say. Now we evaluate these terms.

\subsection{Evaluation of $\I_{1}$}
\label{sect:J1-eval-RH}
Using Lemma~\ref{Laplace-formula} we immediately get 
\begin{equation}
\label{I1-eval-RH}
\I_{1}
 =
\frac{\pi \Gamma(1/k) }{4 k \Gamma(1+1/k)  }  \sum_{n = N+1}^{N + H}  n^{1/k}e^{-n/N} + \Odip{k}{\frac{H}{N}} 
=
\frac{\pi }{4   e} HN^{1/k} +\Odipm{k}{H^2 N^{1/k-1}+ N^{1/k}}.
\end{equation}

\subsection{Estimation of $\I_{5}$}
\label{sect:I5-estim-RH}
Clearly $\I_{5}= \J_5$ of section \ref{sect:J5-estim}. Hence 
we have that 
\begin{equation}
 \label{I5-estim-RH-series} 
 \I_{5}
\ll_k
N^{1/(2k)} ( N^{1/2} + H ) L^2 
\end{equation}
which, comparing with \eqref{I1-eval-RH}, is under control for $H=\infty(N^{1/2-1/(2k)}L^2)$.

\subsection{Estimation of $\I_{4}$}
\label{sect:I4-estim-RH}

Clearly $\I_{4}= \J_4$ of section \ref{sect:J4-estim}. Hence 
we have that 
\begin{equation}
 \label{J4-estim-series-RH}
 \I_{4}
\ll_k 
E(k)
\end{equation}
%
%which, comparing with \eqref{I1-eval-RH}, is under control for $H=\infty(N^{3/4}L)$.
where $E(k)$ is defined in \eqref{Ek-def}.

\subsection{Estimation of $\I_{2}$}
\label{sect:I2-estim-RH}
Now we estimate $\I_{2}$.
Using the identity $f^{2}-g^{2} = 2f(f-g) -(f-g)^{2}$ we obtain 
\begin{equation}
\label{J2-estim1-RH}
\I_{2} 
\ll_k
 \int_{-1/2}^{1/2}
      \vert \Etilde_{2}(\alpha)  \vert
        \frac{\vert U(\alpha,H)\vert}{\vert z\vert^{1/2+1/k}}
     \, \dx \alpha
+
 \int_{-1/2}^{1/2}
       \vert \Etilde_{2}(\alpha)  \vert^{2}  
       \frac{\vert U(\alpha,H) \vert}{\vert z\vert^{1/k}} 
        \, \dx \alpha 
    =
   J_{1}+J_{2},
\end{equation}
say.
Using   \eqref{z-estim},  \eqref{UH-estim}, Lemma \ref{LP-Lemma-gen} and a partial integration argument we obtain
\begin{align}
\notag
J_{2}
&
\ll_k
HN^{1/k}
 \int_{-1/N}^{1/N}  \vert \Etilde_{2}(\alpha)  \vert^{2}   \, \dx \alpha 
 +
 H \int_{1/N}^{1/H}  \vert \Etilde_{2}(\alpha)  \vert^{2}  
 \,\frac{\dx \alpha}{\alpha^{1/k}} 
 +
 \int_{1/H}^{1/2}  \vert \Etilde_{2}(\alpha)  \vert^{2}   
 \,\frac{\dx \alpha}{\alpha^{1+1/k}}
 \\
 \label{J2-estim}
 &\ll_k
HN^{1/k-1/2} L^{2} 
+ 
H^{1/k} N^{1/2} L^{2}  
\ll_k
H^{1/k} N^{1/2} L^{2}  . 
\end{align} 
For $J_1$ we need few cases.
Let $k\ge 3$.
Using  the Cauchy-Schwarz inequality and arguing as for $J_{2}$ we get 
\begin{align}
\notag
J_{1}
&\ll_k
H N^{1/2+1/k} \Bigl(  \int_{-1/N}^{1/N} \dx \alpha  \Bigr)^{1/2}
\Bigl(
 \int_{-1/N}^{1/N} \!\! \vert \Etilde_{2}(\alpha)  \vert^{2}   \, \dx \alpha 
 \Bigr)^{1/2}
 \\ \notag
 & \hskip0.4cm+
 H \Bigl(  \int_{1/N}^{1/H} \!\! \frac{\dx \alpha}{\alpha^{2/k}}\Bigr)^{1/2} 
\Bigl(
 \int_{1/N}^{1/H} \vert \Etilde_{2}(\alpha)  \vert^{2}  
 \frac{\dx \alpha}{\alpha}
  \Bigr)^{1/2} 
     +
  \Bigl(  \int_{1/H}^{1/2}   \frac{\dx \alpha}{\alpha^{2+2/k}}\Bigr)^{1/2}  
\Bigl(
 \int_{1/H}^{1/2}  \vert \Etilde_{2}(\alpha)  \vert^{2}   
 \frac{\dx \alpha}{\alpha}
   \Bigr)^{1/2}
 \\
 \notag
 &\ll_k
H N^{1/k-1/4} L 
+ 
H^{1/2+1/k} N^{1/4} L  \Bigl( 1 + \int_{1/N}^{1/H}  \frac{\dx \xi}{\xi}   \Bigr)^{1/2} 
 +
H^{1/2+1/k} N^{1/4} L  \Bigl( 1+   \int_{1/H}^{1/2}  \frac{\dx \xi}{\xi}   \Bigr)^{1/2}
 \\
 \label{J1-estim-k-RH}
 &\ll_k
H N^{1/k-1/4} L + H^{1/2+1/k} N^{1/4} L^{3/2}
\ll
H^{1/2+1/k} N^{1/4} L^{3/2}. 
\end{align}
For $k=2$ arguing as before we get
\begin{equation}
 \label{J1-estim-2-RH}
J_{1} \ll H N^{1/4} L^{2}.
\end{equation}
Combining \eqref{J2-estim1-RH}-\eqref{J1-estim-2-RH},
and assuming $H\ge N^{1/2}$,
we finally obtain
\begin{equation}
\label{J2-estim-final-RH}
\I_{2} \ll_k H^{1/2+1/k} N^{1/4} L^{2}
\end{equation}
for every $k \ge 2$.

\subsection{Estimation of $\I_{3}$}
\label{sect:I3-estim-RH}
Now we estimate $\I_{3}$.
By the Cauchy-Schwarz  inequality,  \eqref{UH-estim} and Lemma
\ref{Hua-Rieger-Lemma}  we obtain
\begin{align}
\notag
\I_{3}
&\ll_k
\Bigl( \int_{-1/2}^{1/2}  \vert\Stilde_2(\alpha)\vert^{4}   \, \dx \alpha \Bigr)^{1/2}
\Bigl( 
\int_{-1/2}^{1/2}  \vert \Etilde_{k}(\alpha)  \vert^{2} \vert U(\alpha,H)\vert^2 
 \, \dx \alpha 
 \Bigr)^{1/2}
\\
 \label{J3-estim-final-RH}
&\ll_k
N^{1/2}L
\Bigl( H^2
\int_{-1/H}^{1/H}  \vert \Etilde_{k}(\alpha)  \vert^{2} \, \dx \alpha
+ 
 \int_{1/H}^{1/2}  \vert \Etilde_{k}(\alpha)  \vert^{2}   \, \frac{\dx \alpha}{\alpha^{2}} 
 \Bigr)^{1/2}
\ll_k
H^{1/2}N^{1/2+1/(2k)}L^2,
 \end{align}
where in the last step we used  Lemma \ref{LP-Lemma-gen} and 
a partial integration argument.

\subsection{Final words}
By \eqref{main-split-RH}-\eqref{J4-estim-series-RH} and \eqref{J2-estim-final-RH}-\eqref{J3-estim-final-RH},
we can finally write for $H\ge N^{1/2}$ that
\begin{equation}
\label{almost-done-RH}
  \sum_{n = N+1}^{N + H} 
  e^{-n / N} r_k(n)= \frac{\pi }{4   e} HN^{1/k}
+   
\Odipm{k}{H^2 N^{1/k-1}+H^{1/2}N^{1/2+1/(2k)}L^2 + E(k)}
\end{equation}
which is an asymptotic formula for $\infty( N^{1-1/k}L^4) \le H \le \odi{N}$.
{}From  $e^{-n/N}=e^{-1}+ \Odi{H/N}$ for $n\in[N+1,N+H]$, $1\le H \le N$,
we get 
\begin{align*}
   \sum_{n = N+1}^{N + H} 
r_k(n)
&=
 \frac{\pi }{4   e} HN^{1/k}
+
\Odipm{k}{H^2 N^{1/k-1}+H^{1/2}N^{1/2+1/(2k)}L^2 + E(k)}
\\&
+
  \Odipg{k}{\frac{H}{N}\sum_{n = N+1}^{N + H} r_k(n)
}.
\end{align*}
Using $e^{n/N}\le  e^{2}$ 
and \eqref{almost-done-RH},
 the last error term is
$\ll_k H^2N^{1/k-1}+  H^{3/2}N^{-1/2+1/(2k)}L^2 + HN^{-1} E(k)$.
Hence we get
\begin{equation*}
%\label{th-RH-series}
   \sum_{n = N+1}^{N + H} 
r_k(n)
=
 \frac{\pi  }{4   e} HN^{1/k}
+
\Odipm{k}{H^2 N^{1/k-1}+H^{1/2}N^{1/2+1/(2k)}L^2 + E(k)},
\end{equation*}
uniformly for $\infty( N^{1-1/k}L^4) \le H \le \odi{N}$.
Theorem \ref{existence-RH} follows.

%\renewcommand{\bibliofont}{\normalsize}
  
 %%
% \providecommand{\bysame}{\leavevmode\hbox to3em{\hrulefill}\thinspace}
%\providecommand{\MR}{\relax\ifhmode\unskip\space\fi MR }
%% \MRhref is called by the amsart/book/proc definition of \MR.
%\providecommand{\MRhref}[2]{%
%  \href{http://www.ams.org/mathscinet-getitem?mr=#1}{#2}
%}
%\def\ZBL#1{\href{http://zbmath.org/?q=an:#1}{ZBL:#1}} 
%\providecommand{\href}[2]{#2}
%  
%\bibliographystyle{amsplain-nobysame} 
%\bibliography{teonum} 

\begin{thebibliography}{10}

\bibitem{Brudern2016}
J.~Br{\"u}dern, \emph{A ternary problem in additive prime number theory}, in: {From
  Arithmetic to Zeta-Functions. Number Theory in Memory of Wolfgang Schwarz}
  ({J. Sanders {\it et al}}, ed.), Springer, 2016, pp.~57--81.

\bibitem{HarmanK2010}
G.~Harman and A.V. Kumchev, \emph{{On sums of squares of primes II}}, J. Number Theory \textbf{130} (2010), 1969--2002.

\bibitem{Hua1938}
L.~K. Hua, \emph{Some results in the additive prime number theory}, Quart. J.
  Math. Oxford \textbf{9} (1938), 68--80.

\bibitem{LanguascoZ2016b}
A.~Languasco and A.~Zaccagnini, \emph{Short intervals asymptotic formulae for
  binary problems with primes and powers, {II}: density $1$}, Monatsh. Math.
  \textbf{181} (2016), 419--435.

\bibitem{LanguascoZ2016a}
A.~Languasco and A.~Zaccagnini, \emph{Sum of one prime and two squares of
  primes in short intervals}, J. Number Theory \textbf{159} (2016),
  45--58.

\bibitem{LanguascoZ2017b}
A.~Languasco and A.~Zaccagnini,
  \emph{Short intervals asymptotic formulae for binary problems with
  primes and powers, {I}: density $3/2$}, Ramanujan J. \textbf{42} (2017),
  no.~2, 371--383.

\bibitem{LanguascoZ2017c}
A.~Languasco and A.~Zaccagnini,
  \emph{Sums of four prime cubes in short intervals}, Submitted for
  publication. Arxiv preprint 1705.04457, 2017.

\bibitem{LeungL1993}
M.C. Leung and M.C. Liu, \emph{On generalized quadratic equations in three
  prime variables}, Monatsh. Math. \textbf{115} (1993), 133--167.

\bibitem{Li2012}
T.~Li, \emph{On sums of squares of primes and a {$k$}th power of prime}, Rocky
  Mountain J. Math. \textbf{42} (2012), 201--222.

\bibitem{Lu2006}
G.-S. L\"u, \emph{Note on a result of {H}ua}, Adv. Math. (China) \textbf{35}
  (2006), 343--349.

\bibitem{Montgomery1994}
H.~L. Montgomery, \emph{Ten {Lectures} on the {Interface} {Between} {Analytic}
  {Number} {Theory} and {Harmonic} {Analysis}}, CBMS Regional Conference Series
  in Mathematics, vol.~84, A.M.S., 1994.

\bibitem{MontgomeryV1974}
H.~L. Montgomery and R.~C. Vaughan, \emph{Hilbert's inequality}, J. London
  Math. Soc. \textbf{8} (1974), 73--82.
  
\bibitem{Schwarz1961a}
W.~Schwarz, \emph{Zur {D}arstellung von {Z}ahlen durch {S}ummen von
  {P}rimzahlpotenzen. {II}.}, J. Reine Angew. Math. \textbf{206} (1961),
  78--112.

\bibitem{Vaughan1997}
R.~C. Vaughan, \emph{The {Hardy}-{Littlewood} method}, second ed., Cambridge U.
  P., 1997.

\bibitem{Zhao2014a}
L.~Zhao, \emph{{The additive problem with one prime and two squares of
  primes}}, J. Number Theory \textbf{135} (2014), 8--27.

\end{thebibliography}

% 
 
\vskip0.5cm
\noindent
\begin{tabular}{l@{\hskip 20mm}l}
Alessandro Languasco               & Alessandro Zaccagnini\\
Universit\`a di Padova     & Universit\`a di Parma\\
 Dipartimento di Matematica  &  Dipartimento di Scienze Matematiche,\\
 ``Tullio Levi-Civita'' &  Fisiche e Informatiche \\
Via Trieste 63                & Parco Area delle Scienze, 53/a \\
35121 Padova, Italy            & 43124 Parma, Italy\\
{\it e-mail}: alessandro.languasco@unipd.it        & {\it e-mail}:
alessandro.zaccagnini@unipr.it  
\end{tabular}

\end{document}